\documentclass[a4paper,14pt]{article}

\usepackage{multicol}
\usepackage{amsthm}

\usepackage[dvips]{graphicx,color}
\usepackage{amssymb}
\usepackage{amsmath}

\begin{document}

\fontsize{14pt}{16.5pt}\selectfont

\begin{center}
\bf{A dendrite generated from \bf{$\{0,1\}^{\Lambda },~Card\Lambda \succ  \aleph $}}
\end{center}

\fontsize{12pt}{15pt}\selectfont
\begin{center}
Akihiko Kitada$^{1}$, Tomoyuki Yamamoto$^{1,2,*}$, Shousuke Ohmori$^{2}$\\ 
\end{center}

\noindent
$^1$\it{Institute of Condensed-Matter Science, Comprehensive Resaerch Organization, Waseda University,
3-4-1 Okubo, Shinjuku-ku, Tokyo 169-8555, Japan}\\
$^2$\it{Faculty of Science and Engineering, Waseda~University, 3-4-1 Okubo, Shinjuku-ku, Tokyo 169-8555, Japan}\\
~~\\
~~\\

\rm
\noindent
{\bf{Abstract}}\\
\noindent
The existence of a decomposition space with a dendritic structure of a topological space $(\{0,1\}^\Lambda ,\tau_{0}^\Lambda )$ is discussed. Here, $\Lambda $ is any set with the cardinal number $\succ  \aleph , \{0,1\}^{\Lambda }=\{\varphi :\Lambda \rightarrow \{0,1\}\}, \tau_0$ is the discrete topology for $\{0,1\}$ and the topology $\tau_0^{\Lambda }$ for $\{0,1\}^\Lambda $ is the topology with the base $\beta =\{<G_{\lambda _1},\dots,G_{\lambda _n}>~;~G_{\lambda_1}\in \tau_0,\dots,G_{\lambda _n}\in \tau_0, \{\lambda _1,\dots,\lambda _n\}\subset \Lambda ,n\in {\bf N}\}$ where the notation $<E_{\lambda _1},\dots,E_{\lambda _n}>$ concerning the subset $E_{\lambda _i}, i=1,\dots,n$ of $\{0,1\}$ denotes the set $\{\varphi :\Lambda \rightarrow \{0,1\}~;~\varphi (\lambda _1)\in E_{\lambda _1},\dots,\varphi (\lambda _n)\in E_{\lambda _n}, \varphi (\lambda )\in \{0,1\}, \lambda \in \Lambda -\{\lambda _1,\dots,\lambda _n\}\}$.\\
\vspace{-2mm}

\bigskip

\section{Introduction}
~~~~In the present paper, the existence of a dendrite as a decomposition space of a topological space each of whose elements is a map $\varphi $ from a set $\Lambda $ whoes cardinal number is greater than or equal to $\aleph $ to $\{0,1\}$ will be discussed.\\
~~~~Some mathmatical definitions and statements which appear in the text are listed below.\\
A)~[1] Let $(\mathcal{D},\tau(\mathcal{D}))$ be a decomposition 
space of a topological space $(X,\tau)$. That is, $\mathcal{D}$ is\\~~~~~a collection of nonempty, mutually disjoint subsets of $X$ such that $\bigcup \mathcal D = \displaystyle\bigcup _{D\in \mathcal D}D = X$,
\\~~~~and $\tau(\mathcal D)$ is the decomposition topology given by 
$\{\mathcal U\subset \mathcal D ; \bigcup \mathcal U = \displaystyle\bigcup 
_{D\in \mathcal U}D\in \tau\}$. A dec-\\
~~~~omposition $\mathcal D$ of 
$(X,\tau)$ is called usc (upper semi continuous) provided that for 
any el-\\
~~~ ement $D\in \mathcal D$ and for any $U\in \tau$ containing 
$D$, there exists $u\in \tau$ such that $D\subset u$ and\\
~ ~~if $A\cap 
u \not = \phi $ for $A\in \mathcal D$, then $A\subset U$.\\
B)~[2] Let $\mathcal D$ be a decomposition of a set $X$. A subset $Y$ of $X$ which is a union of subcol-\\
~~~~lection of $\mathcal D$, that is, for some $\mathcal U\subset \mathcal D, Y=\bigcup \mathcal U(=\displaystyle\bigcup _{D\in \mathcal U}D)$, is said to be $\mathcal D$-saturated. \\
C)~Let $f:(X,\tau)\rightarrow (Y,\tau')$ be a continuous, closed, onto map. Then, i) [3] The decom-
\\~~~~position space $(\mathcal D_f,\tau(\mathcal D_f))$ of $(X,\tau)$ due to $f$ is usc. Here, the decomposition $\mathcal D_f$ of $X$\\
~~~~is given by $\mathcal D_f=\{f^{-1}(y)\subset X;y\in Y\}$. ii) [4] For any pair $f^{-1}(y)$ and $K\in \Im $ ($K$ is a\\
~~~~closed set of $(X,\tau)$) such that $f^{-1}(y)\cap K=\phi $, there exists a $\mathcal D_f$-saturated $V\in \tau$ sat-\\
~~~~isfying the relations $f^{-1}(y)\subset V$ and $V\cap K=\phi $.\\
D)~[5] Let $X$ be a zero-dimensional (abbreviated as 0-dim), perfect, compact $T_1$-space. \\
~~~~Then, for any compact metric space $Y$, there exists a continuous map $f$ from $X$ onto\\
~~~~$Y$.\\
E)~[6] If $f:(X,\tau)\rightarrow (Y,\tau')$ is a quotient map, that is, an onto map $f$ which is character-\\
~~~~ized by the relation $\tau'=\tau_f(=\{u'\subset Y;f^{-1}(u')\in \tau\})$, then $h:(Y,\tau')\rightarrow (\mathcal D_f,\tau(\mathcal D_f))$,\\
~~~~$y\mapsto f^{-1}(y)$ is a homeomorphism. Namely, the decomposition space $\mathcal D_f$ of $X$ is topolo-\\
~~~~gically equivalent to $Y$.\\
F)~[7] A locally connected, connected, compact metric space containing no simple closed\\
~~~~curve as its subspace is called a dendrite.\\
G)~[8] Let $(X,\tau)$ be compact. If any two distinct points $x$ and $y$ are separated in $X$, that\\
~~~~is, there exists a closed and open (clopen) set $G$ such that $x\in G$ and $y\not \in G$, then $(X,$\\
~~~~$\tau)$ is 0-dim.\\

\section{A dendrite as a decomposition space of a 0-dim, perfect, compact $T_1$-space}
~~~~Let $(X,\tau)$ be a 0-dim, perfect, compact $T_1$-space. Since the dendrite $(Y,\tau_d)$ is a compact metric space by its definition 1,F), according to 1,D), there exists a continuous map $f$ from $(X,\tau )$ onto $(Y,\tau_d)$. Any continuous map from a compact space onto a $T_2$-space is easily verified to be closed. A continuous, closed, onto map is obviously a quotient map. We are convinced from 1,E) that the decomposition space $(\mathcal D_f,\tau(\mathcal D_f))$ of $(X,\tau)$ due to the map $f$ is homeomorphic to the dendrite $(Y,\tau_d)$ under the homeomorphism $h:(Y,\tau_d)\rightarrow (\mathcal D_f,\tau(\mathcal D_f)),y\mapsto f^{-1}(y).$ The decomposition space $(\mathcal D_f,\tau(\mathcal D_f))$ is metrized by means of the homeomorphism $h$ and the metric $d$ of $Y$ as $\rho (D,D')=d(h^{-1}(D),h^{-1}(D')), D,D'\in \mathcal D_f$. That is, the decomposition topology $\tau(\mathcal D_f)$ is identical with the metric topology $\tau_\rho $. Since the metric space $(\mathcal D_f,\tau_\rho )$ is a locally connected, connected, compact metric space which contains no simple closed curve owing to the homeomorphism $h$, it must be a dendrite. 
From 1,C),i), this dendrite $(\mathcal D_f,\tau_\rho )$ based on the continuous, closed, onto map $f$ is usc. Therefore, according to 1,C),ii) each element $f^{-1}(y)$ of this dendrite $\mathcal D_f$ is topologically characterized by the property that, for any pair $f^{-1}(y)$ and $K\in \Im $ such that $f^{-1}(y)\cap K=\phi $, there exists $\mathcal D_f$-saturated $V\in \tau$ such that $f^{-1}(y)\subset V$ and $V\cap K=\phi $. \\
~~\\

\section{Topological space ($\{0,1\}^\Lambda ,\tau_0^\Lambda $)}
~~~~Let $\Lambda $ be a set with $Card\Lambda \succ  \aleph $ and let $\tau_0$ be the discrete topology for $\{0,1\}$. $\tau_0^\Lambda $ is a topology for $\{0,1\}^\Lambda=\{\varphi :\Lambda \rightarrow \{0,1\}\}$ with the base $\beta = \{<G_{\lambda _1},\dots,G_{\lambda _n}> ; G_{\lambda _1}\in \tau_0,\dots,G_{\lambda _n}\in \tau_0,~\{\lambda _1,\dots,\lambda _n\}\subset \Lambda , n\in {\bf N}\}$. Here, the set $\{\varphi :\Lambda \rightarrow \{0,1\} ; \varphi (\lambda _1)\in E_{\lambda _1}, \dots,\varphi (\lambda _n)\in E_{\lambda _n}, \varphi (\lambda )\in \{0,1\}, \lambda \in \Lambda -\{\lambda _{1},\dots,\lambda _{n}\}\}$ concerning the subset $E_{\lambda _i}, i=1,\dots,n$ of $\{0,1\}$, is denoted by $<E_{\lambda _1}, \dots, E_{\lambda _{n}}>$.\\
~~~~In the followings, we will show that the space ($\{0,1\}^\Lambda ,\tau_0^\Lambda $) is 0-dim, perfect, compact $T_1$-space. First, we will check dim$\{0,1\}^\Lambda=0$. Since $(\{0,1\},\tau_0)$ is a compact space, the product space ($\{0,1\}^\Lambda ,\tau_0^\Lambda $) is also compact due to Tychonoff. If $\varphi _p$ and $\varphi _q$ are distinct two points of $\{0,1\}^\Lambda $, there exists an index $\lambda _0\in \Lambda $ such that $\varphi _p(\lambda _0)\not =\varphi _q(\lambda _0)$. Without loss of generality, we can put $\varphi _p(\lambda _0)=0$ and $\varphi _q(\lambda _0)=1$. The open set $<G_{\lambda _0}=\{0\}>=\{\varphi :\Lambda \rightarrow \{0,1\};\varphi (\lambda _0)\in \{0\}, \varphi (\lambda )\in \{0,1\}, \lambda \in \Lambda -\{\lambda _0\}\}$ of $\{0,1\}^\Lambda $ clearly contains the point $\varphi _p$ and does not contain the point $\varphi _q$, that is, 
($\{0,1\}^\Lambda ,\tau_0^\Lambda $) is a $T_1$-space. According to 1,G), it is sufficient for the space ($\{0,1\}^\Lambda ,\tau_0^\Lambda $) to be 0-dim, that the open set $<G_{\lambda _0}=\{0\}>$ is, at the same time, a closed set. Any point $\varphi \in <G_{\lambda _0}=\{0\}>$ is not contained in the open set $<G'_{\lambda _0}=\{1\}>$ and any point $\varphi \not \in <G'_{\lambda _0}=\{1\}>$ must be contained in the set $<G_{\lambda _0}=\{0\}>$. Therefore, the set $<G_{\lambda _0}=\{0\}>$ is equal to the complement of the open set $<G'_{\lambda _0}=\{1\}>$. This means that $ <G_{\lambda _0}=\{0\}>$ is a colsed set. 
Next, to prove that the space ($\{0,1\}^\Lambda ,\tau_0^\Lambda $) is perfect, let us consider $H\in \tau_0^\Lambda -\{\phi \}$. For a point $\varphi _a\in H$, there exists a base $<G_{\lambda _1},\dots,G_{\lambda _n}>$ containing $\varphi _a$ such that $<G_{\lambda _1},\dots,G_{\lambda _n}>\subset H$. Since the set $\Lambda $ is not a finite set, there exists $\bar \lambda\in \Lambda -\{\lambda _1,\dots,\lambda _n\}$ such that $\varphi (\bar \lambda)\in \{0,1\}$. Let $\varphi '_a$ be a point of $\{0,1\}^\Lambda $ defined as $\varphi _a(\lambda )=\varphi' _a(\lambda )$ for $\lambda \not = \bar \lambda$ and $\varphi _a(\bar \lambda)\not =\varphi' _a(\bar \lambda )$. It is clear that $\varphi' _a$ is contained in $<G_{\lambda _1},\dots,G_{\lambda _n}>$ as well as $\varphi _a$. Therefore, the nonempty open set $H$ contains at least two points. This means that the space ($\{0,1\}^\Lambda ,\tau_0^\Lambda $) is perfect.\\

\section{Conclusion}
~~~~Since the space ($\{0,1\}^\Lambda ,\tau_0^\Lambda $) is confirmed to be a 0-dim, perfect, compact $T_1$-space, the discussions in the section 2 can be applied to this space. As a consequence, there exists a decomposition space $\mathcal D$ of ($\{0,1\}^\Lambda ,\tau_0^\Lambda $) which has a dendrite structure. Although the space ($\{0,1\}^\Lambda ,\tau_0^\Lambda $) where $Card\Lambda \succ  \aleph $ is not metrizable[9], its decomposition space $\mathcal D$ is metrizable as a dendite. Finally, we note that there exists a decomposition space of the space ($\{0,1\}^\Lambda ,\tau_0^\Lambda $) one of whose decomposition space has the dendritic structure again (See Appendix). \\
~~\\


\fontsize{12pt}{15pt}\selectfont
\section*{Appendix}

~~~~Using the following lemma, we will confirm the existence of a decomposition space $\mathcal D_f$ of ($\{0,1\}^\Lambda ,\tau_0^\Lambda $) which has a decomposition space with the dendritic structure again.\\

$[$Lemma$]$~~Let $(X,\tau)$ be a 0-dim, perfect, compact $T_0$-space. Then, there exists a non empty 0-dim, perfect, compact $T_0$-proper subspace $(Y,\tau_Y)$ of $(X,\tau)$.
\begin{proof}
Since $X$ is a perfect $T_0$-space, there exist at least two distinct points $x$ and $x'$ and there exists an open set $U$ such that, without loss of generality, $x\in U$ and $x'\not\in U$. Since $X$ is 0-dim, there exists a set $u\in \tau\cap \Im$ which contains $x$ and is contained in $U$. Taking this clopen set $u$ as $Y$, we obtain a 0-dim, perfect, compact proper subset $(Y,\tau_Y)$ of $(X,\tau)$.
\end{proof}
Now for the above space $(X,\tau)$ and its subspace $(Y,\tau_Y)$, let us define a map $f:(X,\tau)\rightarrow (Y,\tau_Y)$ as 
\begin{center}~~~
$f(x)=\left\{
\begin{array}{lcr}
x~,~~x\in Y \\
q~,~~x\in X-Y
\end{array}
\right.$
\end{center}
Here, $q$ is an arbitrarily chosen point of $Y$. Since 0-dim $T_0$-space is necessarily a $T_2$-space, the singleton $\{q\}$ is a closed set. Thus, it is easy to see that the map $f$ is a continuous, closed map from $X$ onto $Y$[10], that is, the map $f$ is a quotient map. Let us note that the decomposition space $\mathcal D_f$ due to this $not~one~to~one$ map $f$ is not the trivial decomposition $\{\{x\};x\in X\}$. From 1,E), $(\mathcal D_f,\tau(\mathcal D_f))$ is homeomorphic to 0-dim, perfect, compact $T_0$-space $(Y,\tau_Y)$, and then, $(\mathcal D_f,\tau(\mathcal D_f))$ must be a 0-dim, perfect, compact $T_0$-space. Since any 0-dim, $T_0$-space is a $T_2$-space, all of the discussions in the text concerning the space ($\{0,1\}^\Lambda ,\tau_0^\Lambda $) can be repeated for the decomposition space $(\mathcal D_f,\tau(\mathcal D_f))$ of ($\{0,1\}^\Lambda ,\tau_0^\Lambda $).\\

\section*{Acknowledgement}
The authors are grateful to Dr. H. Ryo, Dr. Y. Yamashita, Dr. K. Yoshida, Dr. Y. Watayoh and Dr. K. Shikama for useful discussions and encouragements.\\
~~\\
~~\\
\vspace{-4mm}

\section*{References and notes}
\vspace{-7mm}
~~$\\$
[1]~~S.B.Nadler Jr., {\it Continuum theory}, Marcel Dekker 1992, p.36, p38.$\\$
[2]~~Ref.[1], p.38, p.39.$\\$
[3]~~proof of i) It is sufficient for $\mathcal D_f$ to be usc that the natural map $\pi :X\rightarrow \mathcal D_f,x\mapsto D \\$
~~~~~(there exists unique $D\in \mathcal D$ such that $x\in D$) is closed (See, for example Ref.[1],p.39).$\\$
~~~~~Since $f$ is continuous and $f(K)$ is closed for $K\in \Im$, the set $\displaystyle\bigcup _{x\in K}\pi(x)$ must be closed$\\$
~~~~~from the relation $\displaystyle\bigcup _{x\in K}\pi(x)=f^{-1}(f(K))$ which can be easily verified. This means th-$\\$
~~~~~at $\pi$ is a closed map. In fact, $\{\pi (x);x\in A\}$ is closed in $(\mathcal D_f,\tau(\mathcal D_f))\Leftrightarrow \displaystyle\bigcup_{x\in A}\pi(x)\in \Im$.$\\$
[4]~~Ref.[1], p.39.$\\$
[5]~~Ref.[1], p.105, p.106 and p.109 in which the space $X$ is supposed to be a metric space.$\\$
~~~~~But the 0-dim, perfect, compact space $X$ is not necessarily a metric space. It is suff-$\\$
~~~~~icient for the space $X$ to have the property that any compact metric space is a con-$\\$
~~~~~tinuous image of it that $X$ is a $T_1$-space. $\\$
[6]~~Ref.[1], p.46.$\\$
[7]~~a) Ref.[1], p.165, b)~A.Kitada, Y.Ogasawara, T.Yamamoto, Chaos, Solitons \& Frac-\\
~~~~~~tals {\bf 34}, 1732 (2007).$\\$
[8]~~proof of G)~Let $U$ be an open set containing $x$. The complement $U^c$ in $X$ of $U$
 is a$\\$
~~~~~compact set. For each point $y\in U^c$, there exists a clopen set $u(y)$ containing 
$y$ such$\\$
~~~~~that $x\not \in u(y)$. We can choose a finite sub cover $\{u(y_1),\dots,u(y_n)\}$ 
of $U^c$ from the op-$\\$
~~~~~en cover $\{u(y) ; y\in U^c\}$. The clopen set 
$(\displaystyle\bigcup_{i\in \bar n} u(y_i))^c$ containing $x$ is contained in $U$.$\\$
~~~~~$(X,\tau) $ is, therefore, 0-dim.$\\$
[9]~~R.Engelking, {\it General Topology}, Heldermann Verlag Berlin 1989, p.260. If $Card\Lambda \prec \\
~~~~~~\aleph_0 $, the space ($\{0,1\}^\Lambda ,\tau_0^\Lambda $) is well-known to be metrizable. In the present study, it is$\\$
~~~~~emphasized that the metrizability of the space is not necessary for its decomposition$\\$
~~~~~space to be a metric space "dendrite".$\\$
[10]~proof) At the point $x\in Y, f(x)=x$. For any $U\in \tau_Y$ containing $f(x)$, there exists $U\\$
~~~~~containing $x$ such that $f(U)=U\subset U.$ Since $Y$ is an open set, $U\in \tau$. At the point$\\$
~~~~~$x\in Y^c$, for any $U\in \tau_Y$ containing $f(x)=q$, there exists $Y^c\in \tau$ containing $x$ such $\\$
~~~~~that $f(Y^c)=\{q\}\subset U.$ Thus the map $f$ is continuous at any point of $X.$ Next, let $K\\$
~~~~~be a closed set of $X$. $f(K)=f(K\cap X)=f(K\cap (Y\cup Y^c))=f((K\cap Y)\cup (K\cap Y^c))\\
~~~~~~=f(K\cap Y)\cup f(K\cap Y^c)=(K\cap Y)\cup \{q\}\in \Im_Y$. Therefore, $f$ is a closed map.

\end{document}